
  \magnification 1200
  \input amssym  

  \newcount\fontset
  \fontset=1
  \def \dualfont#1#2#3{\font#1=\ifnum\fontset=1 #2\else#3\fi}

  \dualfont\bbfive{bbm5}{cmbx5}
  \dualfont\bbseven{bbm7}{cmbx7}
  \dualfont\bbten{bbm10}{cmbx10}

  \font \eightbf = cmbx8
  \font \eighti = cmmi8 \skewchar \eighti = '177
  \font \eightit = cmti8
  \font \eightrm = cmr8
  \font \eightsl = cmsl8
  \font \eightsy = cmsy8 \skewchar \eightsy = '60
  \font \eighttt = cmtt8 \hyphenchar\eighttt = -1

  \font \sixi = cmmi6 \skewchar \sixi = '177
  \font \sixrm = cmr6
  \font \sixsy = cmsy6 \skewchar \sixsy = '60
  \font \tensc = cmcsc10

  \scriptfont \bffam = \bbseven
  \scriptscriptfont \bffam = \bbfive
  \textfont \bffam = \bbten

  \newskip \ttglue

  \def \eightpoint {\def \rm {\fam0 \eightrm }%
  \textfont0 = \eightrm
  \scriptfont0 = \sixrm \scriptscriptfont0 = \fiverm
  \textfont1 = \eighti
  \scriptfont1 = \sixi \scriptscriptfont1 = \fivei
  \textfont2 = \eightsy
  \scriptfont2 = \sixsy \scriptscriptfont2 = \fivesy
  \textfont3 = \tenex
  \scriptfont3 = \tenex \scriptscriptfont3 = \tenex
  \def \it {\fam \itfam \eightit }%
  \textfont \itfam = \eightit
  \def \sl {\fam \slfam \eightsl }%
  \textfont \slfam = \eightsl
  \def \bf {\fam \bffam \eightbf }%
  \textfont \bffam = \bbseven
  \scriptfont \bffam = \bbfive
  \scriptscriptfont \bffam = \bbfive
  \def \tt {\fam \ttfam \eighttt }%
  \textfont \ttfam = \eighttt
  \tt \ttglue = .5em plus.25em minus.15em
  \normalbaselineskip = 9pt
  \def \MF {{\manual opqr}\-{\manual stuq}}%
  \let \sc = \sixrm
  \let \big = \eightbig
  \setbox \strutbox = \hbox {\vrule height7pt depth2pt width0pt}%
  \normalbaselines \rm }



  \newcount \secno \secno = 0
  \newcount \stno \stno =0
  \newcount \eqcntr \eqcntr=0

  \def \ifn #1{\expandafter \ifx \csname #1\endcsname \relax }

  \def \track #1#2#3{\ifn{#1}\else {\tt\ [#2 \string #3] }\fi}

  \def \laberr#1#2{\message{*** RELABEL CHECKED FALSE for #1 ***}
      RELABEL CHECKED FALSE FOR #1, EXITING.
      \end}

  \def \seqnumbering {\global \advance \stno by 1 \global
    \eqcntr=0 \number \secno .\number \stno }

  \def \current {\number \secno
    \ifnum \number \stno = 0\else .\number \stno \fi }

  \def \eqmark #1 {\global \advance\eqcntr by 1
    \edef\a{\current.\number\eqcntr}
    \eqno {(\a)}
    \syslabel{#1}{\a}
    \track{showlabel}{*}{#1}}

  \def \syslabel#1#2{\global \expandafter \edef \csname
    #1\endcsname {#2}}

  \def \fcite#1#2{\syslabel{#1}{#2}\lcite{#2}}

  \def \label #1 {%
    \ifn {#1}%
      \syslabel{#1}{\current}%
    \else
      \edef\a{\expandafter\csname #1\endcsname}%
      \edef\b{\current}%
      \ifx \a \b \else \laberr{#1=(\a)=(\b)} \fi
      \fi
    \track{showlabel}{*}{#1}}

  \def \lcite #1{(#1\track{showcit}{$\bullet$}{#1})}

  \def \cite #1{[{\bf #1}\track{showref}{\#}{#1}]}

  \def \scite #1#2{{\rm [\bf #1\track{showref}{\#}{#1}{\rm \hskip 0.7pt:\hskip 2pt #2}\rm]}}


 \def \Headlines #1#2{\nopagenumbers
    \advance \voffset by 2\baselineskip
    \advance \vsize by -\voffset
    \headline {\ifnum \pageno = 1 \hfil
    \else \ifodd \pageno \tensc \hfil \lcase {#1} \hfil \folio
    \else \tensc \folio \hfil \lcase {#2} \hfil
    \fi \fi }}

  \def \Date #1 {\footnote {}{\eightit Date: #1.}}


  \def \lcase #1{\edef \auxvar {\lowercase {#1}}\auxvar }

  \def \goodbreak {\vskip0pt plus.1\vsize \penalty -250 \vskip0pt
plus-.1\vsize }

  \def \section #1{\global\def \SectionName{#1}\stno = 0 \global
\advance \secno by 1 \bigskip \bigskip \goodbreak \noindent {\bf
\number \secno .\enspace #1.}\medskip \noindent \ignorespaces}

  \long \def \sysstate #1#2#3{\medbreak \noindent {\bf \seqnumbering
.\enspace #1.\enspace }{#2#3\vskip 0pt}\medbreak }
  \def \state #1 #2\par {\sysstate {#1}{\sl }{#2}}
  \def \definition #1\par {\sysstate {Definition}{\rm }{#1}}
  \def \remark #1\par {\sysstate {Remark}{\rm }{#1}}


  \def \proof {\medbreak \noindent {\it Proof.\enspace }}
  \def \proofend {\ifmmode \eqno \square \else \hfill \square
\looseness = -1 \medbreak \fi }

  \def \$#1{#1 $$$$ #1}
  \def \=#1{\buildrel \hbox{\sixrm #1} \over =}

  \def \Item #1{\smallskip \item {#1}}
  \newcount \zitemno \zitemno = 0

  \def \izitem {\zitemno = 0}
  \def \zitemplus {\global \advance \zitemno by 1\relax}
  \def \rzitem{\romannumeral \zitemno}
  \def \rzitemplus {\zitemplus \rzitem}
  \def \zitem {\Item {{\rm(\rzitemplus)}}}

  \newcount \nitemno \nitemno = 0
  
  \def \nitem {\global \advance \nitemno by 1 \Item {{\rm(\number\nitemno)}}}

  \newcount \aitemno \aitemno = 0
  \def\boxlet#1{\hbox to 6.5pt{\hfill #1\hfill}}
  \def \iaitem {\aitemno = 0}
  \def \aitem {\Item {(\ifcase \aitemno \boxlet a\or \boxlet b\or
\boxlet c\or \boxlet d\or \boxlet e\or \boxlet f\or \boxlet g\or
\boxlet h\or \boxlet i\or \boxlet j\or \boxlet k\or \boxlet l\or
\boxlet m\or \boxlet n\or \boxlet o\or \boxlet p\or \boxlet q\or
\boxlet r\or \boxlet s\or \boxlet t\or \boxlet u\or \boxlet v\or
\boxlet w\or \boxlet x\or \boxlet y\or \boxlet z\else zzz\fi)} \global
\advance \aitemno by 1}

  \newcount \footno \footno = 1
  \newcount \halffootno \footno = 1
  \def \footcntr {\global \advance \footno by 1
  \halffootno =\footno
  \divide \halffootno by 2
  $^{\number\halffootno}$}


  \def \N {{\bf N}}
  
  \def \<{\left \langle \vrule width 0pt depth 0pt height 8pt }
  \def \>{\right \rangle }

  \def \and {\hbox {,\quad and \quad }}

  \def \imply {\mathrel{\Rightarrow}}
  \def \for #1{,\quad \forall\,#1}
  \def \square {\hbox {$\sqcap \!\!\!\!\sqcup $}}
  
  \def \stress #1{{\it #1}\/}
  
  \def \*{\otimes}

  \newcount \bibno \bibno =0
  \def \newbib #1{\global\advance\bibno by 1 \edef #1{\number\bibno}}
  \def \bibitem #1#2#3#4{\smallskip \item {[#1]} #2, ``#3'', #4.}
  \def \references {
    \begingroup
    \bigskip \bigskip \goodbreak
    \eightpoint
    \centerline {\tensc References}
    \nobreak \medskip \frenchspacing }

  \def \proofend {\ifmmode \eqno \square \else \hfill $\square$
\looseness = -1 \medbreak \fi }

  \def\emptyset{\varnothing}

  \font\rs=rsfs10
  
  \def\leq{\leqslant}
  \def\geq{\geqslant}

  \font\mf=cmex10
  \def\union {\mathop{\raise 9pt \hbox{\mf S}}\limits}
  \def\inters{\mathop{\raise 9pt \hbox{\mf T}}\limits}
  \def\bigcup{\union}
  \def\bigcap{\inters}
  \def\imply{\ \Rightarrow \ }
  \def\dil{\mathrel{|}}
  \def\phi{\varphi}

  \def\S{\hbox {\rs G}\kern2.5pt} \def\S{\Lambda}
  \def\St{\S^{(2)}}
  \def\Su{\tilde{\S}}
  \def\D#1{\S^{#1}}

  \def \reqmark #1 {\eqno {(\seqnumbering)}
    \syslabel{#1}{\current}
    \track{showlabel}{*}{#1}}
  \def\case#1{\medskip\noindent {\tensc Case #1}:}
  \def\its{\Cap}
  \def\disj{\perp}
  \def\A#1#2{A(#1,#2)}

  \def\rep#1{S_{#1}}  
  \def\fin#1{P_{#1}}
  \def\ini#1{Q_{#1}}
  \def\d{s}
  \def\r{r}
  \def\obj#1{{\rm obj}(#1)}
  \def\Gen{{\cal G}}
  \def\a{\alpha}
  \def\b{\beta}
  \def\g{\gamma}
  \def\CS{{\cal O}(\S)}
  \def\CSU{\tilde\CS}
  \def\HHGA{C^*(\S)}


  \def \newbib #1#2{\global\advance\bibno by 1 \edef
#1{\number\bibno}}

\newbib\BHRS{BHRS}
\newbib\BPRS{BPRS}
\newbib\Blackadar{Blackadar}
\newbib\inverse{inverse}
\newbib\actions{actions}
\newbib\infinoa{infinoa}
\newbib\Muhly{Muhly}
\newbib\FLR{FLR}
\newbib\Katsura{Katsura}
\newbib\KPActions{KPActions}
\newbib\KP{KP}
\newbib\KPR{KPR}
\newbib\KPRR{KPRR}
\newbib\PQR{PQR}
\newbib\PRRS{PRRS}
\newbib\PatGraph{PatGraph}
\newbib\Renault{Renault}
\newbib\RaeBook{RaeBook}
\newbib\RSY{RSY}
\newbib\RaeSzy{RaeSzy}
\newbib\PTW{PTW}
\newbib\RSa{RSa}
\newbib\RSb{RSb}
\newbib\Tomforde{Tomforde}
\newbib\Watatani{Watatani}
\newbib\Yeend{Yeend}

  \Headlines
  {Semigroupoid C*-Algebras}
  {R.~Exel}

  \null\vskip -1cm
  \centerline{\bf SEMIGROUPOID C*-ALGEBRAS}
  \footnote{\null}
  {\eightrm 2000 \eightsl Mathematics Subject Classification:
  \eightrm 
  Primary 46L05; 
  secondary 18B40. 
  }

  \bigskip
  \centerline{\tensc 
    R.~Exel\footnote{*}{\eightpoint Partially supported by
CNPq.}}

  \bigskip
  \Date{29 Nov 2006 (revised 15 Mar 2007)}

  \midinsert 
  \narrower \narrower
  \eightpoint \noindent
  A semigroupoid is a set equipped with a partially defined
associative operation.  Given a semigroupoid $\S$ we construct a
C*-algebra $\CS$ from it.  We then present two main examples of
semigroupoids, namely the Markov semigroupoid associated to an
infinite 0--1 matrix, and the semigroupoid associated to a row-finite
higher-rank graph without sources.  In both cases the semigroupoid
C*-algebra is shown to be isomorphic to the algebras usually attached
to the corresponding combinatorial object, namely the Cuntz-Krieger
algebras and the higher-rank graph C*-algebras, respectively.
  In the case of a higher-rank graph $(\S,d)$, it follows that the
dimension function $d$ is superfluous for defining the corresponding
C*-algebra.  
  \endinsert

  \section{Introduction}
  The theory of C*-algebras has greatly benefited from Combinatorics
in the sense that some of the most interesting examples of C*-algebras
arise from combinatorial objects, such as the case of graph
C*-algebras 
  \cite{%
  \BHRS,
  \BPRS,
  \FLR,
  \KPActions,
  \KPR,
  \KPRR, 
  \PatGraph,
  \RaeSzy,
  \PTW,
  \Tomforde, 
  \Watatani
  },  
  see also \cite{\RaeBook} and the references therein.  More
recently Kumjian and Pask have introduced the notion of
\stress{higher-rank graphs} \cite{\KP}, inspired by Robertson and
Steger's work on buildings \cite{\RSa}, \cite{\RSb}, which turns out
to be another combinatorial object with which an interesting new class
of C*-algebras may be constructed.  See also 
  \cite{%
  \Muhly, 
  \PQR, 
  \PRRS, 
  \RaeBook, 
  \RSY}.

The crucial insight leading to the notion of higher-rank graphs lies
in viewing ordinary graphs as \stress{categories} (in which the
morphisms are finite paths) equipped with a \stress{length} function
possessing a certain unique factorization property (see \cite{\KP} for
more details).

Kumjian and Pask's interesting idea of viewing graphs as categories
suggests that one could construct C*-algebras for more general
categories.  

Since Eilenberg and Mac Lane introduced the notion of
categories in the early 40's, the archetypal idea of composition of
functions has been mathematically formulated in terms of 
categories, whereby a great emphasis is put on the \stress{domain} and
\stress{co-domain} of a function.  However one may argue that, while
the domain is an intrinsic part of a function, co-domains are
not so vital.  If one imagines a very elementary function $f$ with
domain, say $X=\{1,2\}$, defined by $f(x) = x^2$, one does not really
need to worry about its co-domain.  But if $f$ is to be seen as a
morphism in the category of sets, one needs to first choose a set $Y$
containing the image of $f$, and only then $f$ becomes an element of
Hom$(X,Y)$.  Regardless of the very innocent nature of our
function $f$, it suddenly is made to evoke an enormous amount of
morphisms, all of them having the same domain $X$, but with the
wildest possible collection of co-domains.  

Addressing this concern
one could replace the idea of categories with the following: a big set
(or perhaps a class) would represent the collection of all morphisms,
regardless of domains, ranges or co-domains.  A set of
\stress{composable} pairs $(f,g)$ of morphisms would be given in
advance and for each such pair one would define a composition $fg$.
Assuming the appropriate associativity axiom one arrives as the notion
of a \stress{semigroupoid}, precisely defined in
\fcite{DefineSemiGroupoid}{2.1} below.  Should our morphisms be actual
functions a sensible condition for a pair $(f,g)$ to be composable
would be to require the range of $g$ to be contained in the domain of
$f$,  but we might also think of more abstract situations in which the
morphisms are not necessarily functions.  

For example, let $A = \{\A ij\}_{i,j\in\Gen}$, be an infinite 0-1
matrix, where $\Gen$ is an arbitrary set, and let $\S_A$ be the set of
all finite \stress{admissible words} in $\Gen$, meaning finite sequences
  $
  \a = \a_1\a_2\ldots \a_n,
  $
  of elements $\a_i\in\Gen$, such that $\A{\a_i}{\a_{i+1}}=1$.
  Given $\a, \b\in\S_A$  write 
  $$
  \a = \a_1\a_2\ldots \a_n
  \and
  \b = \b_1\b_2\ldots \b_m,
  $$
  and let us say that $\a$ and $\b$ are composable if
$\A{\a_n}{\b_1}=1$, in which case we let $\a\b$ be the
concatenated word
  $$
  \a\b = \a_1\a_2\ldots \a_n \b_1\b_2\ldots \b_m,
  $$
  We shall refer to $\S_A$ as the \stress{Markov} semigroupoid.  
  This category-like structure lacks a notion of \stress{objects} and
in fact it cannot always be made into a category.  Consider for
instance the matrix
  $$
  A=\pmatrix{1 & 1\cr 1 & 0}.
  $$
  If we let the index set of $A$ be $\Gen=\{\a_1,\a_2\}$, notice that
the words $\a_1$ and $\a_2$ may be legally composed to form the words
  $\a_1\a_1$, $\a_1\a_2$, and $\a_2\a_1$, but $\a_2\a_2$ is forbidden,
precisely because  $\A {\a_2}{\a_2}=0$.

Should there exist an underlying category, the fact that, say,
$\a_1\a_2$ is a legal composition would lead one to believe that
$\d(\a_1)$, the domain, or source of $\a_1$ coincides with $\r(\a_2)$,
the co-domain of $\a_2$.
  But then for similar reasons one would have
  $$
  \d(\a_2) = \r(\a_1) = \d(\a_1) = \r(\a_2),
  $$
  which would imply that $\a_2\a_2$ is a valid composition, but it is
clearly not.  This example was in fact already noticed by Tomforde
\cite{\Tomforde} with the purpose of showing that a 0-1 matrix $A$ is
not always the edge matrix of a graph.
  Although the above matrix may be replaced by another one which is the
edge matrix of a graph and gives the same Cuntz-Krieger algebra, the
same trick does not work for infinite matrices.

This is perhaps an indication that we should learn to live with
semigroupoids which are not true categories.
  Given the sheer simplicity of the notion of semigroupoid, one can
easily fit all of the combinatorial objects so far referred to within
the framework of semigroupoids.

The goal of this work is therefore to introduce a notion of
\stress{representation} of semigroupoids, with its accompanying
universal C*-algebra, which in turn generalizes earlier constructions
such as the Cuntz-Krieger algebras for arbitrary matrices of
\cite{\infinoa} and the higher-rank graph C*-algebras of \cite{\KP},
and hence ordinary graph C*-algebras as well.

The present paper is mostly devoted to comparing semigroupoid
C*-algebras with Cuntz-Krieger and higher-rank graph C*-algebras.
Please see \cite{\actions}, where a deeper study is made of the
structure of semigroupoid C*-algebras, including describing them
as groupoid C*-algebras.

The definition of a representation of a semigroupoid $\S$ given in
\fcite{DefineRepTwo}{4.1}, and consequently of the C*-algebra of $\S$,
here denoted $\CS$, is strongly influenced by \cite{\infinoa}, and
hence it is capable of smoothly dealing with the troubles usually
caused by \stress{non-row-finiteness}.

  Speaking of another phenomenon that requires special attention in
graph C*-algebra theory, the presence of \stress{sources}, once cast
in the perspective of semigroupoids, becomes much easier to deal with.

To avoid confusion we use a different term and define a
\stress{spring} (rather than source) to be an element $f$ of a
semigroupoid $\S$ for which $fg$ is not a legal multiplication for any
$g\in\S$.  The sources of graph theory are much the same as our
springs, and they cause the same sort of problems, but there are some
subtle, albeit important differences.  For example, in a category any
element $f$ may be right-multiplied by the identity morphism on its
domain, and hence categories never have any springs.  On the other
hand, even though higher-rank graphs are defined as categories,
sources may still be present and require a special treatment.
See however \scite{\actions}{18.2.ii}.

While springs are irremediably killed when considered within the
associated semigroupoid C*-algebra, as shown in
\fcite{SourcesVanish}{5.1}, it is rather easy to get rid of them by
replacing the given semigroupoid by a somewhat canonical spring-less
one \fcite{SpringToSource}{3.3}.
  This is specially interesting because a slight correction performed
on the ingredient semigroupoid is seen to avoid the need to redesign
the whole theoretical apparatus.

As already mentioned, the C*-algebra $\HHGA$ associated to a
higher-rank graph $(\S,d)$ in \cite{\KP} turns out to be a special
case of our construction: since $\S$ is defined to be a category, it
is obviously a semigroupoid, so we may consider its semigroupoid
C*-algebra $\CS$, which we prove to be isomorphic to $\HHGA$ in
\fcite{HHGvsSGPG}{8.7}.

One of the most interesting aspects of this is that the construction
of $\CS$ does not use the \stress{dimension function} ``$d$\kern 1.5pt" at all,
relying  exclusively on the algebraic structure of the subjacent
category.  In other words, this shows that the dimension function is
superfluous in the definition of $\HHGA$.

It should be stressed that our proof of the isomorphism between $\CS$
and $\HHGA$ is done under the standing hypotheses of \cite{\KP},
namely that $(\S,d)$ is row-finite and has no sources.  The reader
will not find here a comparison between our construction and the more
recent treatment of Farthing, Muhly and Yeend \cite{\Muhly} for
general finitely aligned higher rank graphs.  We hope to be able to
address this issue in a future paper.

It is a consequence of  Definition \lcite{\DefineRepTwo},
describing our notion of a  representation
$\rep\null$ of a semigroupoid $\S$,  that if
$\S$ contains elements $f$, $g$ and $h$ such that $fg=fh$, then 
  $$
  \rep g = \rep h.
  $$
  Therefore, even if $g$ and $h$ are different, that difference is
blurred when these  elements are seen in $\CS$ via the
universal representation.
  This should probably be interpreted as saying that our
representation theory is not really well suited to deal with general
semigroupoids in which non monic elements are present.  
An element $f$ is said to be \stress{monic} if
  $$
  fg=fh \imply g=h.
  $$
  Fortunately
all of our examples consist of semigroupoids containing only monic
elements.
  See section (4) for more details.

No attempt has been made to consider \stress{topological}
semigroupoids although we believe this is a worthwhile program to be
pursued.  Among a few indications that this can be done is Katsura's
topological graphs \cite{\Katsura} and Yeend's \cite{\Yeend}
topological higher-rank graphs, not to mention Renault's pioneering
work on groupoids \cite{\Renault}.

After recognizing the precise obstruction for interpreting
Cuntz-Krieger algebras from the point of view of categories or graphs,
one can hardly help but to think of the obvious generalization of
higher-rank graphs to semigroupoids based on the unique factorization
property.  Even though we do not do anything useful based on this
concept we spell out the precise definition in \fcite{HHSgpd}{8.1}
below.  As an example, the Markov semigroupoid for the above $2\times
2$ matrix is a rank 1 semigroupoid which is not a rank 1 graph.

We would also like to mention that although we have not seriously
considered the ultra-graph C*-algebras of Tomforde \cite{\Tomforde}
from a semigroupoid point of view, we
believe that these may also be described in terms of naturally
occurring semigroupoids.

I would like to acknowledge many fruitful
conversations with
  A.~Kumjian,
  M.~Laca and   
  D.~Pask
  during the process of developing this work.
  Special thanks go to A.~Sims for bringing to our attention some
important references in the subject of higher-rank graphs.

\section{Semigroupoids}
  In this section we introduce the basic algebraic ingredient of our
construction.
  
\definition
  \label DefineSemiGroupoid
  A semigroupoid is a triple $(\S,\St,\ \cdot\ )$ such that $\S$ is a set,
$\St$ is a subset of $\S\times\S$, and
  $$
  \cdot\ : \St \to \S
  $$
  is an operation which is associative in the following sense:  if
$f,g,h\in\S$ are such that either
  \begingroup \parindent 25pt
  \izitem 
  \zitem $(f,g)\in\St$ and $(g,h)\in\St$, or \global\edef\BasicForm{\rzitem}
  \zitem $(f,g)\in\St$ and $(f g,h)\in\St$, or \global\edef\LeftForm{\rzitem}
  \zitem $(g,h)\in\St$ and $(f,\ g h)\in\St$,

  \medskip\noindent 
  then all of $(f,g)$, $(g,h)$, $(f g,h)$ and $(f,g h)$ lie in
$\St$, and
  $$
  (f g) h = f (g h).
  $$
  Moreover, for every $f\in\S$, we will let
  $$
  \D f = \left\{g\in\S: (f,g)\in\St\right\}.
  $$
  \endgroup

From now on we fix a semigroupoid $\S$.

\definition 
  \label DefineDivision
  Let $f,g\in\S$. We shall say that $f$ \stress{divides}
$g$, or that $g$ is a \stress{multiple} of $f$, in symbols $f\dil g$,
if either 
  \izitem
  \zitem $f=g$, or
  \zitem  there exists $h\in\S$ such that
  $
  fh=g.
  $
  \medskip\noindent
  When $f\dil g$, and $g\dil f$, we shall say that $f$ and $g$ are
\stress{equivalent}, in symbols $f\simeq g$.

Perhaps the correct way to write up the above definition is to require
that $(f,h)\in\St$ before referring to the product ``$fh$".  However
we will adopt the convention that, when a statement is  made
about a freshly introduced element which  involves a
multiplication, then the statement is implicitly supposed to include
the requirement that the multiplication involved is allowed.

Notice that in the absence of anything resembling a unit in $\S$, it
is conceivable that for some element $f\in\S$ there exists no $u\in\S$
such that $f=fu$.  Had we not explicitly included
\lcite{\DefineDivision.i}, it would not always be  the case that $f\dil f$.

A useful artifice is to introduce a unit for $\S$, that is, pick some
element in the universe outside $\S$, call it $1$, and set
$\Su=\S\mathop{\dot\cup}\ \{1\}$.  For every $f\in\S$ put
  $$ 
  1f = f1 = f.
  $$
  Then, whenever $f\dil g$, regardless of whether $f=g$ or not, there
always exists $x\in\Su$ such that $g=fx$.
 
We will find it useful to extend the definition of $\D f$, for
$f\in\Su$, by putting
  $$
  \D 1 = \S.
  $$
  Nevertheless, even if $f1$ is a meaningful product for every
$f\in\S$, we \stress{will not} include $1$ in $\D f$.

We should be aware that $\Su$ is \stress{not} a semigroupoid.  Otherwise, 
since $f1$ and $1g$ are meaningful products, axiom
\lcite{\DefineSemiGroupoid.\BasicForm} would imply that $(f1)g$ is
also a meaningful product, but this is clearly not always the case.

It is interesting to understand the extent to which the associativity
property fails for $\Su$.  As already observed,
\lcite{\DefineSemiGroupoid.i} does fail irremediably when $g=1$.
Nevertheless it is easy to see that \lcite{\DefineSemiGroupoid}
generalizes to $\Su$ in all other cases.  This is quite useful, since
when we are developping a computation, having arrived at an expression
of the form $(fg)h$, and therefore having already checked that all
products involved are meaningful, we most often want to proceed by
writing
  $$
  \ldots = (fg)h = f(gh).
  $$
  The axiom to be invoked here is \lcite{\DefineSemiGroupoid.ii} (or
\lcite{\DefineSemiGroupoid.iii} in a similar situation), and
fortunately not \lcite{\DefineSemiGroupoid.i}!

\state Proposition
  Division is a reflexive and transitive relation.

\proof
   That division is reflexive follows from the definition. In order to
prove transitivity let $f,g\in\S$ be such that $f\dil g$ and $g\dil
h$.  We must prove that $f\dil h$.

 The case in which $f=g$, or $g=h$ is obvious.  Otherwise 
there are $u,v$ in $\S$ (rather than in $\Su$) such that $fu=g$, and
$gv=h$.  As observed above, it is implicit that $(f,u),(g,v)\in\St$,
which implies that
  $$
  (f,u),(fu,v)\in\St.
  $$
  By \lcite{\DefineSemiGroupoid.\LeftForm} we deduce that $(u,v)\in\St$ and
that
  $$
  f(uv) = (fu) v = gv = h,
  $$
  and hence $f\dil h$.  
  \proofend

Division is also invariant under multiplication on the left:

\state Proposition
  \label DivisionLeftInvar
  If $k,f,g\in\S$ are such that 
  $f\dil g$, and
  $(k,f)\in\St$, then 
  $(k,g)\in\St$ and $kf\dil kg$.

  \proof The case in which $f=g$ being obvious we 
assume that there is $u\in\S$  such that $fu=g$. 
  Since $(k,f),\ (f,u)\in\St$  we conclude from
\lcite{\DefineSemiGroupoid.\BasicForm} that 
  $(kf,u)$ and $ (k,g)=(k,fu)$
  lie in $\St$, and that
  $$
  (kf) u = k (fu) = kg,
  $$
  so $kf\dil kg$.
  \proofend

The next concept will be crucial to the analysis of the structure of
semigroupoids.

\definition 
  Let $f,g\in\S$.  We shall say that $f$ and $g$ \stress{intersect} if
they admit a \stress{common multiple}, that is, an element $m\in\S$
such that $f\dil m$ and $g\dil m$.  Otherwise we will say that $f$ and
$g$ are \stress{disjoint}.  We shall indicate the fact that $f$ and
$g$ intersect by writing $f\its g$, and when they are disjoint we will
write $f\disj g$.

If there exists  a right-zero element, that is, an element $0\in\S$ such
that $(f,0)\in\St$ and $f0=0$, for all $f\in\S$, then obviously $f\dil 0$,
and hence any two elements  intersect.  We shall be mostly interested
in semigroupoids without a right-zero element.

Employing  the unitization $\Su$ notice that $f\its g$ if and only if
there are  $x,y\in\Su$ such that $fx=gy$.

A last important concept, borrowed from the Theory of Categories, is
as follows:

\definition
  We shall say that an element $f\in\S$ is \stress{monic} if for every
$g,h\in\S$ we have
  $$
  fg=fh \imply g=h.
  $$

  \def\Ga{\Gamma}
  \def\Gat{\Ga^{(2)}}

\section{Springs}
  We would now like to discuss certain special properties of elements
$f\in\S$ for which $\D f=\emptyset$.  It would be sensible to call
these elements \stress{sources}, 
  following the terminology adopted in Graph Theory, but given some
subtle differences we'd rather use  another term:

\definition \label DefineSpring
  We will say that an element $f$ of a semigroupoid $\S$ is a
\stress{spring} when $\D f=\emptyset$.

Springs  are sometimes annoying, so we shall now discuss a way of
getting rid of springs.  Let us therefore fix a semigroupoid $\S$
which has springs. 

Denote by $\S_0$ the subset of $\S$ formed by all springs and let $E'$
be a set containing a distinct element $e'_g$, for every $g\in\S_0$.
Consider any equivalence relation ``$\sim$" on $E'$ according to which
  $$
  e'_g\sim e'_{fg},
  \eqno{(\seqnumbering)}
  \label SpringInvariance
  $$
  for any spring $g$, and any $f$ such that $g\in\D f$.  Observe that
$fg$ is necessarily also a spring since $\D {fg} = \D g$, by
\lcite{\DefineSemiGroupoid.i-ii}.
  For example, one can take the equivalence relation according to which
any two elements are related.  Alternatively we could use
the smallest equivalence relation satisfying
\lcite{\SpringInvariance}.

We shall denote the quotient space \ $E'/\kern-3pt\sim$ \ by $E$, and
for every spring $g$ we will denote the equivalence class of $e'_g$
by $e_g$.
  Unlike the $e'_g$, the $e_g$ are obviously no longer distinct elements.
In particular we have 
  $$
  e_g =e_{fg}
  \for f\in\S
  \for g\in\S_0\cap \D f.
  $$
  We shall now construct a semigroupoid $\Ga$ as follows:  set 
  $
  \Ga = \S \mathop{\dot\cup} E,
  $
  and put 
  $$
  \Gat = \St \ \cup \ \big\{(g,e_g): g\in\S_0\big\} \ \cup \
\big\{(e_g,e_g): g\in\S_0\big\}.
  $$
  Define the multiplication 
  $$
  \cdot\ : \Gat \to \Ga,
  $$
  to coincide with the multiplication of $\S$ when restricted to
$\St$, and moreover set
  $$
  g\cdot e_g = g 
  \and 
  e_g\cdot e_g=e_g 
  \for g\in\S_0.
  $$
  It is rather tedious, but entirely elementary, to show that $\Ga$ is
a semigroupoid without any springs containing $\S$.  
  To summarize the conclusions of this section we state the following:

\state Theorem 
  \label SpringToSpring
  For any semigroupoid $\S$ there exists a spring-less semigroupoid\/
$\Ga$ containing $\S$.

Given a certain freedom in the choice of the equivalence relation
``$\sim$" above, there seems not to be a canonical way to embed $\S$
in a spring-less semigroupoid.  The user might therefore have to make a
case by case choice according to his or her preference.

\section{Representations of semigroupoids}
  \label RepresenSect
  In this section we begin the study of the central notion bridging
semigroupoids and operator algebras.

\definition
  \label DefineRepTwo
  Let  $\S$  be a semigroupoid and let $B$ be a unital C*-algebra.
  A mapping 
  $
  \rep\null : \S \to B
  $
  will be called a \stress{representation of $\S$ in $B$}, if for
every $f,g\in\S$, one has that:
  \izitem 
  \zitem $\rep f$ is a partial isometry,
  \zitem
  $
  \rep f\rep g  = \left\{\matrix{
    \rep {fg}, &\hbox{ if } (f,g)\in\St, \cr
    \vrule height 15pt width 0pt
    0, & \hbox{ otherwise.}\hfill}\right.
  $
  \medskip\noindent
  Moreover the initial projections $\ini f = \rep f^*\rep f$, and the
final projections $\fin g=\rep g\rep g^*$, are required to commute
amongst themselves and to satisfy
  \medskip
  \zitem $\fin f\fin g=0$, if $f\disj g$, 
  \zitem $\ini f\fin g = \fin g$, if $(f,g)\in\St$. 

\medskip Notice that if  $(f,g)\notin\St$, then
  $
  \ini f\fin g = \rep f^*\rep f \rep g\rep g^* = 0,
  $
  by (ii).  Complementing (iv) above we could therefore add:

\medskip
 \zitem $\ini f\fin g = 0$, if $(f,g)\notin\St$.

\bigskip
  We will automatically extend any representation $\rep\null$ to the
unitization $\Su$ by setting $\rep 1 = 1$.  Likewise we put $\ini
1=\fin 1=1$.

Notice that in case $\S$ contains an element $f$ which is not monic,
say $fg=fh$, for a pair of distinct elements $g,h\in\S$, one
necessarily has $\rep g=\rep h$, for every representation
$\rep\null$.  In fact
  $$
  \rep g =
  \rep g\rep g^*\rep g =
  \fin  g\rep g = 
  \ini f \fin  g\rep g = 
  \rep f^*\rep f   \rep g\rep g^*\rep g =
  \rep f^*\rep f   \rep g =
  \rep f^*\rep {fg},
  $$
  and similarly $\rep h = \rep f^*\rep {fh}$, so it follows that $\rep
g=\rep h$, as claimed.
 
This should probably be interpreted as saying that our representation
theory is not really well suited  to deal with general semigroupoids
in which non monic elements are present.  In fact, all of our examples
consist of semigroupoids containing only monic elements.

  From now on we will fix a representation $\rep\null$ of a given
semigroupoid $\S$ in a unital C*-algebra $B$.
  By \lcite{\DefineRepTwo.iv} we have that
  $
  \fin h \leq \ini f,
  $
  for all  $h\in \D f$,
  so if $h_1,h_2\in \D f$ we deduce that
  $$
  P_{h_1} \vee P_{h_2} := P_{h_1} + P_{h_2} - P_{h_1}P_{h_2} \leq \ini f.
  $$
  More generally, 
  if $H$ is a finite subset of $\D f$ we will have
  $$
  \bigvee_{h\in H}  \fin h \leq \ini f.
  $$
  We now wish to discuss whether or not the above inequality becomes
an identity under circumstances  which we now make explicit:

\definition
  Let $X$ be any subset of $\S$.  A subset $H\subseteq X$ will be
called a \stress{covering} of $X$ if for every $f\in X$ there exists
$h\in H$ such that $h\its f$.  If moreover the elements of $H$ are
mutually disjoint then $H$ will be called a \stress{partition} of $X$.

The following elementary fact is noted for further reference:

  \state Proposition
  \label PartitionsAreMaximal
  A subset $H\subseteq X$ is a partition of $X$ if and only if $H$ is
a maximal subset of $X$ consisting  of pairwise  disjoint elements.

Returning to our discussion above we wish to require that 
  $$
  \bigvee_{h\in H} \fin h \=?  \ini f,
  \eqno{(\seqnumbering)}
  \label FirstWish
  $$
  whenever $H$ is a covering of $\D f$.
  The trouble with this equation is that when $H$ is infinite there is
no reasonable topology available on $B$ under which one can make sense
of the supremum of infinitely many commuting projections.

Before we try to attach any sense to \lcite{\FirstWish} notice that if
$g\in\Su$ and $h\in \S\setminus \D g$, then $\fin h\leq 1- \ini g$, by
\lcite{\DefineRepTwo.v}, and hence also
  $$
  \bigvee_{h\in H}\fin h \leq 1-\ini g,
  $$
  for every finite set $H\subseteq \S\setminus\D g$.
  More generally, 
  given finite subsets $F,G\subseteq\Su$,  denote by 
  $$
  \S^{F,G} =
  \Big(\bigcap_{f\in F} \D f\Big) \cap
  \Big(\bigcap_{g\in G} \S\setminus\D g\Big),
  $$
  and let $h\in \S^{F,G}$.  By \lcite{\DefineRepTwo.iv-v}, we have
that
  $$
  \fin h \leq
  \prod_{f\in F}\ini f\prod_{g\in G}(1-\ini g). 
  $$
  As in the above cases we deduce that 
  $$
  \bigvee_{h\in H}\fin h \leq
  \prod_{f\in F}\ini f\prod_{g\in G}(1-\ini g),
  $$
  for every finite subset  $H\subseteq\S^{F,G}$.

\definition
  \label DefineTight
  A representation $\rep\null$ of $\S$ in a unital C*-algebra $B$ is
said to be \stress{tight} if for every finite subsets
$F,G\subseteq\Su$, and for every finite covering $H$ of
$\S^{F,G}$ one has that
  $$
  \bigvee_{h\in H}\fin h =
  \prod_{f\in F}\ini f \prod_{g\in G}(1-\ini g). 
  $$
 
Observe that if no such covering exists, then any representation is
tight by default.

For almost every representation theory there is a C*-algebra whose
representations are in one-to-one correspondence with the
representations in the given theory.  Semigroupoid representations are no
exception:

\definition
  \label DefineCStarAlgSGPD
  Given a semigroupoid $\S$ we shall let $\CSU$ be the universal
unital C*-algebra generated by a family of partial isometries $\{\rep
f\}_{f\in\S}$ subject to the relations that the correspondence
$f\mapsto\rep f$ is a tight representation of $\S$.  That
representation will be called the \stress{universal representation}
and the closed *-subalgebra of $\CSU$ generated by its range will be
denoted $\CS$.

It is clear that $\CSU$ is either equal to $\CS$ or to its
unitization.
 Observe also that the relations we are referring to in the above
definition are all expressable in the form described in
\cite{\Blackadar}.  Moreover these relations are admissible, since
any partial isometry has norm one.  It therefore follows that $\CSU$
exists.

The universal property of $\CSU$ may be expressed as follows:

\state Proposition
  For every tight representation $T$ of $\S$ in a unital
C*-algebra $B$ there exists a unique *-homomorphism
  $$
  \phi:\CSU  \to B,
  $$
 such that $\phi(\rep f)=T_f$, for every $f\in\S$.

It might also be interesting to define a ``Toeplitz" extension of
$\CSU$, as the universal unital C*-algebra generated by a family of
partial isometries $\{\rep f\}_{f\in\S}$ subject to the relations that
the correspondence $f\mapsto\rep f$ is a (not necessarily tight)
representation of $\S$. If such an algebra is denoted ${\cal T}(\S)$,
it is immediate that $\CSU$ is a quotient of ${\cal T}(\S)$.

As already observed the usefulness of these constructions is probably
limited to the case in which every element of $\S$ is monic.

\section{Tight representations and springs}
  Tight representations and springs do not go together well, as
explained below:

\state Proposition 
  \label  SpringsVanish
  Let $\rep\null$ be a tight representation of a semigroupoid $\S$ and
let $f\in\S$ be a spring (as defined in \lcite{\DefineSpring}).  Then
$\rep f=0$.

  \proof
  Under the assumption that $\D f=\emptyset$, notice that the empty set
is a covering of $\D f$ and hence $\ini f=0$, by \lcite{\DefineTight}.  Since
$\ini f = \rep f^*\rep f$, one has that $\rep f=0$, as well.
  \proofend

We thus see that springs do not play any role with respect to tight
representations.  There are in fact some other non-spring elements on
which every tight representation vanishes.  Consider for instance the
situation in which $\D f$ consists of a finite number of elements, say
$\D f=\{h_1,\ldots,h_n\}$, each $h_i$ being a spring.  Then $\D
f$ is a finite cover of itself and hence by \lcite{\DefineTight} we
have
  $$
  \ini f = \bigvee_{i=1}^n \fin {h_i} = 0,
  $$
  which clearly implies that $\rep f=0$.

One might feel tempted to redesign the whole concept of tight
representations especially if one is bothered by the fact that springs
are killed by them.  However we strongly feel that the right thing to
do is to redesign the semigroupoid instead, using
\lcite{\SpringToSpring} to replace $\S$ by a spring-less semigroupoid
containing it.

In this case it might be useful to understand the following
situation:

\state Proposition
  Let $\rep\null$ be a tight representation of a semigroupoid $\S$ and
suppose that $f\in\S$ is such that $\D f$ contains a single element
$e$ such that $e^2=e$.  Then  $\rep e$ is a projection and moreover  $\rep e =
\ini f$.

\proof
  Since $e^2=e$, we have that $\rep e^2=\rep e$.  But since $\rep e$
is also a partial isometry, it must necessarily be a projection.
  By assumption we have that $\{e\}$ is a finite covering for $\D f$
so
  $$
  \ini f = \fin e = \rep e \rep e^* = \rep e.
  \proofend
  $$

With this in mind we will occasionally work under the assumption that
our semigroupoid has no springs.

\def\TCKCond#1{T\kern-0.8ptC\kern-0.5ptK$_{#1}$}
\def\OAu{\tilde {\cal O}_{\kern-2pt A}}
\def\SA{\S}
\def\CSA{\CSU}
\def\SASmall{\SA}
\def\SAt{\SA^{(2)}}
\def\SAXY{\SA^{X,Y}}
\def\SAu{\tilde\SA}
\def\DA#1{\SA^{#1}}

\def\accenta#1{\check{#1}}
\def\accentc#1{\hat{#1}}
\def\repa#1{\accenta{S}_{#1}}
\def\repc#1{\accentc{S}_{#1}}  

\def\inia#1{\accenta{Q}_{#1}}
\def\inic#1{\accentc{Q}_{#1}}
\def\fina#1{\accenta{P}_{#1}}
\def\finc#1{\accentc{P}_{#1}}

\section{The Markov semigroupoid}
  In this section we shall present a semigroupoid whose C*-algebra is
isomorphic to the Cuntz-Krieger algebra introduced in \cite{\infinoa}.
  For this let $\Gen$ be any set and let $A = \{\A ij\}_{i,j\in \Gen}$
be an arbitrary matrix with entries in $\{0,1\}$.
  We consider the set $\SA=\SA_A$ of all finite \stress{admissible}
words
  $$
  \a = \a_1\a_2\ldots \a_n,
  $$
  i.e., finite sequences of elements $\a_i\in\Gen$, such that
$\A{\a_i}{\a_{i+1}}=1$.  Even though it is sometimes interesting to
consider the empty word as valid, we shall \stress{not} do so.  If
allowed, the empty word would duplicate the role of the extra element
$1\in\SAu$. Our words are therefore assumed to have strictly positive
length ($n\geq1$).

 Given another admissible word, say
  $\b = \b_1\b_2\ldots \b_m$,
  the concatenated word 
  $$
  \a\b:= \a_1\ldots \a_n \b_1\ldots \b_m
  $$
  is admissible as long as  $\A{\a_n}{\b_1}=1$.  Thus, if we set
  $$
  \SAt=\{(\a,\b)=(\a_1\a_2\ldots \a_n,\ \b_1\b_2\ldots
\b_m)\in\SA\times\SA: \A{\a_n}{\b_1}=1\},
  $$ 
  we get a semigroupoid with concatenation as product.  

  \definition
  The semigroupoid $\SA=\SA_A$ defined above will be called the
\stress{Markov} semigroupoid.

Observe that the springs in $\SA$ are precisely the words
$\a=(\a_1\a_2\ldots \a_n)$ for which $A(\a_n,j)=0$, for every
$j\in\Gen$, that is, for which the $\a_n^{th}$ row of $A$ is zero.
  To avoid springs we will assume that no row of $A$ is zero.

\state Theorem
  \label IsoOACA
  Suppose that $A$ has no zero rows. Then $\CSA$ is *-isomorphic to
the unital Cuntz-Krieger algebra $\OAu$ defined in
\scite{\infinoa}{7.1}.

  \proof Throughout this proof we will denote the standard generators
of $\OAu$ by $\{\repa x\}_{x\in\Gen}$, together with their initial and
final projections
  $\inia x = \repa x^*\repa x$ and 
  $\fina x= \repa x\repa x^*$, respectively.
  Likewise the standard generators of $\CSA$ will be denoted by
$\{\repc f\}_{f\in\SASmall}$, along with with their initial and final
projections
  $\inic f = \repc f^*\repc f$ and 
  $\finc f= \repc f\repc f^*$.
  In addition, for every $x\in\Gen$ we will identify the one-letter
word ``$x$" with the element $x$ itself, so we may think of $\Gen$ as
a subset of $\SA$.

  We begin by claiming that the set of partial isometries
  $$
  \{\repc x\}_{x\in\Gen}\subseteq \CSA
  $$
  satisfies the defining relations of $\OAu$, namely \TCKCond1,
\TCKCond2, and \TCKCond3 of \scite{\infinoa}{Section 3}, plus
\scite{\infinoa}{1.3}.

Conditions \TCKCond1 and \TCKCond2 follow immediately from 
\lcite{\DefineRepTwo}, and the observation that if $x$ and $y$ are
distinct elements of $\Gen$, then $x\disj y$ as elements of $\SA$.

When $\A ij=1$ we have  that $(i,j)\in\SAt$ and hence 
  $
  \finc i\inic j = \finc j = \A ij \finc j,
  $
  by \lcite{\DefineRepTwo.iv}.
  Otherwise, if $\A ij=0$, we have that $(i,j)\notin\SAt$ and
hence
  $
  \finc i\inic j = 0= \A ij \finc j,
  $
  by  \lcite{\DefineRepTwo.v}.  This proves
\TCKCond3.

In order to prove \scite{\infinoa}{1.3} let $X,Y$ be finite subsets of
$\Gen$ such that
  $$
  A(X,Y,j) := \prod_{x\in X}\A xj \prod_{y\in Y} (1-\A yj)
  \eqmark AXYj
  $$
  equals zero for all but finitely many  $j$'s.  It is then easy to
see that 
  $$
  Z := \{j\in\Gen : A(X,Y,j)\neq 0 \}
  $$
  is a finite partition of $\SAXY$, so 
  $$
  \prod_{x\in X}\inic x \prod_{y\in Y}(1-\inic y) = 
  \bigvee_{j\in Z}\finc j = 
  \sum_{j\in Z}\finc j = 
  \sum_{j\in \Gen}A(X,Y,j)\finc j,
  $$
  because the canonical representation $f\in\SA \mapsto \repc f\in
\CSA$ is tight by definition.
   It then follows from the universal property of $\OAu$ that there
exists a *-homomorphism
  $$
  \Phi : \OAu \to \CSA,
  \eqmark OneHomomorphism
  $$
  such that $\Phi(\repa x) = \repc x$, for every $x\in\Gen$.  

  Next consider the map
  $
  \repa\null:\SA \to \OAu
  $
  defined as follows:  given $\a\in\SA$, write   $\a = \a_1\a_2\ldots
\a_n$, with $\a_i\in\Gen$, and put 
  $$
  \repa 
  \a = \repa{\a_1}\repa{\a_2}\ldots\repa{\a_n}.
  $$
  We claim that $\repa\null$ is a tight representation of $\SA$ in
$\OAu$.  The first two axioms of \lcite{\DefineRepTwo} are immediate,
while the commutativity of the $\fina f$, and $\inia g$ follow from
\scite{\infinoa}{3.2} and \scite{\inverse}{2.4.iii}.
  Next suppose that $\a,\b\in\SA$ are such that $\a\disj\b$.  One may
then prove that
  $$
  \a=\a_1\ldots\a_p\ldots \a_n
  \and
  \b=\b_1\ldots\b_p\ldots \b_m,
  $$
  with $1\leq p\leq n,m$, and such that $\a_i=\b_i$ for $i<p$, and
$\a_p\neq\b_p$.  Denoting by $\g=\a_1\ldots\a_{p-1}$ (possibly the empty
word), we have that
  $$
  \fina \a \leq \repa\g \repa{\a_p} \repa{\a_p}^* \repa\g^* = 
  \repa\g \fina{\a_p} \repa\g^*,
  $$
  and similarly
  $ 
  \fina \b \leq  \repa\g \fina{\b_p} \repa\g^*.
  $
  It follows that
  $$
  \fina \a\fina\b \leq
  \repa\g \fina{\a_p} \repa\g^* \repa\g \fina{\b_p} \repa\g^* =
  \repa\g \fina{\a_p} \inia\g \fina{\b_p} \repa\g^* =
  \repa\g \fina{\a_p} \fina{\b_p} \inia\g  \repa\g^* =
  0,
  $$
  by \scite{\infinoa}{\TCKCond2}, hence proving
\lcite{\DefineRepTwo.iii}.  In order to verify
\lcite{\DefineRepTwo.iv} let $(\a,\b)\in\SAt$, so that
$\A{\a_n}{\b_1}=1$, where $n$ is the length of $\a$. 
  As shown in ``Claim 1" in the proof of \scite{\infinoa}{3.2}, we
have that
  $
  \inia \a = \inia{\a_n},
  $
  so
  $$
  \inia \a \fina \b =
  \inia{\a_n} \fina{\b_1}\fina \b =
  \fina{\b_1}\fina \b =
  \fina \b,
  $$
  where we have used \TCKCond3 in the second equality. 

We are then left with the task of proving $\repa\null$ to be tight.
For this let $X$ and $Y$ be finite subsets of $\SA$ and let $Z$ be
a finite covering of
  $
  \SAXY.
  $
  We must prove that
  $$
  \bigvee_{h\in Z}\fina h =
  \prod_{f\in X}\inia f \prod_{g\in Y}(1-\inia g). 
  \eqmark EquationToBeProved
  $$
  Using \TCKCond3 it is easy to check the inequality ``$\leq$" in
\lcite{\EquationToBeProved} so it suffices to verify the opposite
inequality.

Let $h_1,h_2\in Z$ be such that $h_1\its h_2$, and write
$h_1x_1=h_2x_2$, where $x_1,x_2\in\SAu$.  Assuming that the length of
$h_1$ does not exceed that of $h_2$, one sees that $h_1$ is an initial
segment of $h_2$, and hence $h_1\dil h_2$.  Any element of $\SAXY$
which intersects $h_2$ must therefore also intersect $h_1$.  This said
we see that $Z':= Z\setminus\{h_2\}$ is also a covering of
  $
  \SAXY.
  $
  Since the left-hand-side of \lcite{\EquationToBeProved} decreases
upon replacing $Z$ by $Z'$, it is clearly enough to prove the
remaining inequality ``$\geq$" with $Z'$ in place of $Z$.

Proceeding in such a way every time we find pairs of intersecting elements in $Z$
we may then suppose that $Z$ consists of pairwise disjoint elements,
and hence that $Z$ is a partition.

  Given $f\in \SA$, write $f=\a_1\ldots\a_n$, with $\a_i\in\Gen$, and
observe that
  $\inia f = \inia {\a_n}$, as already mentioned.  Since $\DA f =
\DA{\a_n}$, as well, we may assume without loss of generality that $X$
and $Y$ consist of words of length one, or equivalently that
$X,Y\subseteq\Gen$.
  Let 
  $$
  J = \{j\in\Gen : A(X,Y,j)\neq 0 \},
  $$
  where $A(X,Y,j)$ is as in \lcite{\AXYj}.
  Notice that $j\in J$ \ if and only if $\A xj=1$, and $\A yj=0$, for
all $x\in X$ and $y\in Y$, which is precisely to say that $j\in\SAXY$.
In other words
  $$
  J  =  \SAXY\cap \Gen.
  $$ 
  It is clear that $J$ shares with $Z$ the property of being maximal
among the subsets of pairwise disjoint elements of $\SAXY$ (see
\lcite{\PartitionsAreMaximal}).

  Suppose for the moment that $Z$ is formed by words of length one,
i.e, that $Z\subseteq\Gen$.  Then $Z\subseteq J$, and so $Z=J$, by
maximality.  This implies that $J$ is finite and
  $$
  \bigvee_{z\in Z}\fina z =
  \bigvee_{j\in J}\fina j =
  \sum_{j\in J}\fina j =
  \sum_{j\in \Gen}A(X,Y,j)\fina j =
  \prod_{x\in X}\inia x \prod_{y\in Y}(1-\inia y),
  $$
  by \scite{\infinoa}{1.3}, thus proving \lcite{\EquationToBeProved}.
  Addressing the situation in which $Z$ is not necessarily contained in
$\Gen$, let
  $$
  Z_j = \{\a\in Z: \a_1 = j\}
  \for j\in J.
  $$
  Since $Z\subseteq \SAXY$ it is evident that 
  $$
  Z= \bigcup_{j\in J} Z_j.
  $$
  Moreover notice that each $Z_j$ is nonempty since otherwise
$Z\cup\{j\}$ will be a subset of $\SAXY$ formed by mutually disjoint
elements, contradicting the maximality of $Z$.
  In particular this shows that $J$ is finite and hence we may use
\scite{\infinoa}{1.3}, so that
  $$
  \prod_{x\in X}\inia x \prod_{y\in Y}(1-\inia y) = 
  \sum_{j\in \Gen}A(X,Y,j)\fina j =
  \sum_{j\in J}\fina j.
  \eqmark FaltaDescerNaArvore
  $$

We claim that for every $j\in J$ one has that
  $$
  \fina j =
  \sum_{z\in Z_j}  \fina z.
  $$
  Before proving the claim lets us notice that it does implies our
goal, for then 
  $$
  \sum_{z\in Z}  \fina z =
  \sum_{j\in J}\ \sum_{z\in Z_j}  \fina z =
  \sum_{j\in J} \fina j \={(\FaltaDescerNaArvore)}
  \prod_{x\in X}\inia x \prod_{y\in Y}(1-\inia y),
  $$
  proving \lcite{\EquationToBeProved}.
  
Noticing that each $Z_j$ is maximal among subsets of mutually disjoint
elements beginning in $j$,  the claim follows from the following:

\state Lemma
  Given  $x\in\Gen$, let 
  $
  \SA(x) = \{\a\in\SA: \a_1 = x\},
  $
  and let $H$ be a finite partition of $\SA(x)$.  Then
  $$
  \sum_{h\in H} \fina h = \fina x.
  $$
  
  \proof
  Let $n$ be the maximum length of the elements of $H$.  We will prove
the statement by induction on $n$.  If $n=1$ it is clear that  $H=\{x\}$ and the
conclusion follows by obvious reasons.  Supposing that $n>1$ observe
that $x\notin H$, or else  any element in $H$ with length $n$ will
intersect $x$, violating the hypothesis that $H$ consists of mutually
disjoint elements.  Therefore every element of $H$ has length at least two.
 
  Let $J=\{j\in\Gen: \A xj=1\}$ and set 
  $H_j=\{\a\in H: \a_2 = j\}$.  It is clear that 
  $$
  H=\bigcup_{j\in J} H_j.
  $$
  Moreover notice that every $H_j$ is nonempty, since otherwise 
  $H\cup\{xj\}$ consists of mutually disjoint elements and properly
contains $H$, contradicting maximality.   In particular this implies
that $J$ is finite and hence 
by \scite{\infinoa}{1.3} we have
  $$
  \inia x = \sum_{j\in\Gen}\A xj \fina j = \sum_{j\in J} \fina j.
  \eqmark FoundIni
  $$
  For every $j\in J$, let $H_j'$ be the set obtained by deleting the
first letter from all words in $H_j$, so that $H_j'\subseteq \SA(j)$,
and $H_j = xH_j'$.  One moment of reflexion will convince the reader
that $H_j'$ is maximal among the subsets of mutually disjoint elements
of $\SA(j)$.  Since the maximum length of elements in $H_j'$ is no
bigger than $n-1$, we may use induction to conclude that
  $$
  \fina j = \sum_{k\in H_j'} \fina k.
  $$
  Therefore 
  $$
  \fin x = 
  \repa x \repa x^*   \repa x \repa x^* =
  \repa x \inia x  \repa x^* \={(\FoundIni)}
  \sum_{j\in J}  \repa x \fina j \repa x^* =
  \sum_{j\in J}\ \sum_{k\in H_j'} \repa x  \fina k \repa x^* \$=
  \sum_{j\in J}\ \sum_{k\in H_j'} \fina {xk} =
  \sum_{j\in J}\ \sum_{h\in H_j} \fina h =
  \sum_{h\in H} \fina h.
  \proofend
  $$
 
Returning to the proof of \lcite{\IsoOACA}, now in possession  of the
information that $\repa\null$ is a tight representation of $\SA$,
we conclude by the universal property of $\CSA$ that there exists
a *-homomorphism
  $$
  \Psi : \CSA \to \OAu,
  $$
  such that $\Psi(\repc \a)=\repa \a$, for all $\a\in\SA$.  It is then
clear that $\Psi$ is the inverse of the homomorphism $\Phi$ of
\lcite{\OneHomomorphism}, and hence both $\Phi$ and $\Psi$ are
isomorphisms.
  \proofend

\section{Categories}
  \label CategorySection
  In this section we fix a small category $\S$.  Notice that the
collection of all morphisms of $\S$ (which we identify with $\S$
itself) is a semigroupoid under composition.  We shall now study $\S$
from the point of view of the theory introduced in the previous sections.

Given $v\in\obj\S$ (meaning the set of objects of $\S$) we will
identify $v$ with the identity morphism on $v$, so that we will see
$\obj\S$ as a subset of the set of all morphisms.

Given $f\in\S$ we will denote by $\d(f)$ and $\r(f)$ the
\stress{domain} and \stress{co-domain} of $f$, respectively.
  Thus the set of all composable pairs may be described as
  $$
  \St  = \{(f,g)\in\S\times\S: \d(f)=\r(g)\}.
  $$
  Given  $f\in\S$ notice that 
  $
  \D f = \{g\in\S: \d(f)=\r(g)\}.
  $
  In particular, if $v\in \obj\S$ then   $\d(v)=\r(v)=v$, so 
  $$
  \D v = \{g\in\S: \r(g)=v\}.
  $$

A category is a special sort of semigroupoid in several ways.  For
example, if $f_1$, $f_2$, $g_1$ and $g_2\in\S$ are such that
$(f_i,g_i)\in\St$ for all $i,j$, except perhaps for $(i,j)=(2,2)$,
then necessarily $(f_2,g_2)\in\St$, because
  $$
  \d(f_2) = \r(g_1) = \d(f_1) = \r(g_2).
  $$

Another special property of a category among semigroupoids is the fact
that for every $f\in\S$ there exists $u\in\D f$ such that $f = fu$,
namely one may take $u$ to be (the identity on) $\d(f)$.  Thus $f\dil
f$ even if we had omitted \lcite{\DefineDivision.i} in the definition
of division.  Clearly this also implies that $\S$ has no springs.

From now on we fix a representation $\rep\null$ of $\S$ in a unital
C*-algebra $B$ and denote by $\ini f$ and $\fin f$, the initial
and final projections of each $\rep f$, respectively.
  A few elementary facts are in order:

\state Proposition
  \label ElmntFactsOnRepCat
  \izitem 
  \zitem For every $v\in\obj\S$ one has that $\rep v$ is a projection,
and hence $\rep v = \fin v = \ini v$.
  \zitem If $u$ and $v$ are distinct objects then $\fin u\perp\fin v$.
  \zitem  For every  $f\in\S$ one has that  $\ini f = \fin {\d(f)}$.

  \proof
  We leave the elementary proof of (i) to the reader.
  Given distinct objects
  $u$ and $v$ it is clear that $u\disj v$, so $\fin u\perp\fin v$, by
  \lcite{\DefineRepTwo.iii}.
  With respect to (iii) 
  we have
  $$
  \ini f =
  \rep f^*\rep f =
  \rep f^*\rep {f\d(f)} =
  \rep f^*\rep f\rep {\d(f)} =
  \ini f\fin {\d(f)} =
  \fin {\d(f)},
  $$
  where the last equality follows from  \lcite{\DefineRepTwo.iv}.
  \proofend

\definition
  \label DefineNonDegen
  Let $H$ be a Hilbert space and let $\rep\null:\S\to {\cal B}(H)$ be
a representation.  We will say that $\S$ is \stress{nondegenerated}
if the closed *-subalgebra of ${\cal B}(H)$ generated by the range of
$\rep\null$ is nondegenerated.

Nondegenerated Hilbert space representations are partly tight in
the following sense:

\state  Proposition
  \label Partlytight
  Let $\rep\null:\S\to {\cal B}(H)$ be a representation.  If either
  \izitem
  \zitem $\rep\null$ is nondegenerated, or 
  \zitem $\obj\S$ is infinite, 
  \medskip \noindent then for every finite subsets $F,G\subseteq\S$
such that
  $
  \S^{F,G}  = \emptyset,
  $
  one has that
  $$
  \prod_{f\in F}\ini f \prod_{g\in G}(1-\ini g) = 0. 
  $$

  \proof
  Notice that $\S^{F,G}=\emptyset$,  implies that
  $$
  \Big(\bigcap_{f\in F} \D f\Big) \subseteq
  \S \setminus   \Big(\bigcap_{g\in G} \S\setminus\D g\Big) =
  \bigcup_{g\in G} \D g.
  \eqmark CompleComplementar
  $$

\medskip \noindent{\tensc Case 1:} Assuming that $F\neq\emptyset$,
let $f_0\in F$.  Then either there is some $f\in F$, with
$\d(f)\neq\d(f_0)$, in which case
  $$
  \ini {f_0} \ini f = \fin {\d(f_0)} \fin {\d(f)} = 0,
  $$
  proving the statement; or $\d(f)=\d(f_0)$, for all $f\in F$.
Therefore we may suppose that $\d(f_0)$ belongs to $\D f$ for every
$f\in F$, and hence by \lcite{\CompleComplementar} there exists
$g_0\in G$ such that $\d(f_0)\in\D {g_0}$.  But this is only possible
if $\d(f_0)=\d(g_0)$ and hence
  $$
  \ini {f_0} (1 - \ini{g_0}) =
  \fin{\d(f_0)} (1-\fin{\d(g_0)}) = 0,
  $$  
  concluding the proof in  case 1.

\medskip \noindent{\tensc Case 2:} Assuming next that $F=\emptyset$, we
claim that
  $$
  \obj\S=\{\d(g): g\in G\}.
  $$
  In fact,   arguing as in \lcite{\CompleComplementar} one has that 
  $
  \bigcup_{g\in G} \D g = \S,
  $
  so for every $v\in\obj\S$ there exists $g$ in $G$ such that
$v\in\D g$, whence  $v=\d(g)$, proving our claim. 

Under the assumption that $\obj\S$ is infinite we have reached a
contradiction, meaning that case 2 is impossible and the proof is
concluded.  We thus proceed supposing nondegeneracy.  Let
  $$
  R = \prod_{g\in G}(1-\ini g),
  $$
  so, proving the statement is equivalent to proving that $R=0$.
  Given $v\in\obj\S$, let $g\in G$ be such that $v=\d(g)$.  Then
  $$
  \ini g \rep v =
  \fin {\d(g)} \rep v =
  \rep v,
  $$
  from where we deduce that
  $$
  R\rep v = 
  R(1-\ini g) \rep v = 0.
  $$
  Given any $f\in\S$ we then have that
  $$
  R\rep f = R\rep {\r(f)} \rep f = 0
  \and
  R\rep f^* =
  R \rep {\d(f)}  \rep f^* = 0,
  $$
  so $R=0$, by nondegeneracy.
  \proofend

We next present a greatly simplified way to check that a
representation of $\S$ is tight.

\state  Proposition
  \label tightInCategories
  Given a representation $\rep\null:\S\to {\cal B}(H)$, consider the
following two statements:
  \iaitem
  \aitem $\rep\null$ is tight.
  \aitem For every $v\in\obj\S$ and every finite covering $H$ of\/ $\D
v$ one has that
  $
  \bigvee_{h\in H} \fin h = \fin v.
  $
  \medskip\noindent Then
  \izitem
  \zitem (a) implies (b). 
  \zitem If $\rep\null$ is nondegenerated, or $\obj\S$ is infinite,
then (b) implies (a).

  \proof
  (i): Assume that $\rep\null$ is tight and that $H$ is a finite covering of
$\D v$.  Setting $F=\{v\}$ and $G=\emptyset$, notice that
  $$
  \S^{F,G} = \D v,
  $$
  so $H$ is a finite covering of $\S^{F,G}$, and hence we have by
definition that
  $$
  \bigvee_{h\in H} \fin h = 
  \prod_{f\in F}\ini f\prod_{g\in G}(1-\ini g) =
  \ini v = \fin v.
  $$

\noindent (ii): 
Assuming $\rep\null$ nondegenerated, or $\obj\S$ infinite, we next
prove that (b) implies (a).
So  let $F$ and $G$ be finite subsets of $\S$ and
let $H$ be a finite covering of
  $
  \S^{F,G}.
  $
  We must prove that the identity in \lcite{\DefineTight} holds.
  If $\S^{F,G}=\emptyset$, the conclusion follows from
\lcite{\Partlytight}.  So we assume that $\S^{F,G}\neq\emptyset$.

\case 1 $F\neq\emptyset$.   
 Pick $h\in \S^{F,G}$ and
notice that for every $f\in F$ one has that $\d(f)=\r(h)$, and for
every $g\in G$, it is the case that $\d(g)\neq\r(h)$.  It therefore
follows that 
  $$
  \S^{F,G} = \D v,  
  $$
  where $v=\r(h)$, so $H$ is in fact a covering of $\D v$.  By
hypothesis we then have that
  $$
  \bigvee_{h\in H} \fin h = \fin v.
  \eqmark SumForPv
  $$
  On the other hand observe that for every $g\in G$, we have that
  $$
  \ini g \={(\ElmntFactsOnRepCat.iii)} \fin {\d(g)} \perp \fin v,
  $$
  given that $\d(g)\neq v$.  Noticing that for 
$f\in F$, we have 
  $
  \ini f  =\fin {\d(f)} = \fin v,
  $
  we deduce that
  $$
  \prod_{f\in F}\ini f\prod_{g\in G}(1-\ini g) = \fin v
\={(\SumForPv)} 
  \bigvee_{h\in H} \fin h,
  $$ 
  proving that  the identity in \lcite{\DefineTight}  indeed holds in
case $F\neq\emptyset$.

\case 2 $F=\emptyset$.  Let
  $$
  V=\obj\S\setminus \{\d(g): g\in G\}, 
  $$
  so that
  $$
  \S^{F,G} = \union_{v\in V} \D v.
  $$
  Given that $H$ is a finite covering of $\S^{F,G}$, we have that for
each $v\in V$ there exists $h\in H$ such that $v\its h$, which in turn
implies that $\r(h)=v$.  Therefore $V$ is finite and hence so is 
$\obj\S$.

Thus, case 2 is impossible under the hypothesis that $\obj\S$ is
infinite, and hence the proof is finished under that hypothesis.
  We therefore proceed supposing nondegeneracy.
It is then  easy to show that 
  $$
  \sum_{v\in \obj\S} \fin v =1,
  $$
  and hence 
  $$
  \prod_{g\in G}(1-\ini g) \={(\ElmntFactsOnRepCat.iii)} 
  \prod_{g\in G}(1-\fin {\d(g)}) =
  \sum_{v\in V} \fin v.
  \eqmark PartialComputation
  $$
  By assumption $H$ is contained in $\S^{F,G}$, and hence the range of
each $h\in H$ belongs to $V$.  Thus
  $$
  H= \dot{\union_{v\in V}} H_v,
  $$
  where $H_v = \{h\in H: \r(h)=v\}$.  Observe that $H_v$ is a
covering for $\D v$, since if $f\in\D v$, there exists some $h\in H$
with $h\its f$, but this implies that $\r(h)=\r(f)=v$, and hence $h\in
H_v$.  Thus
  $$
  \bigvee_{h\in H} \fin h =
  \bigvee_{v\in V}\Big(\bigvee_{h\in H_v} \fin h\Big) =
  \bigvee_{v\in V}\fin v =
  \sum_{v\in V} \fin v \={(\PartialComputation)}
  \prod_{g\in G}(1-\ini g).
  \proofend
  $$

  \def\d{s}
  \def\r{r}

  \def\accenta#1{\check{#1}}
  \def\accentc#1{\hat{#1}}
  \def\repa#1{\accenta{S}_{#1}}
  \def\repc#1{\accentc{S}_{#1}}

  \def\inia#1{\accenta{Q}_{#1}}
  \def\inic#1{\accentc{Q}_{#1}}
  \def\fina#1{\accenta{P}_{#1}}
  \def\finc#1{\accentc{P}_{#1}}

\section{Higher-rank graphs}
  We shall now apply the conclusions above to show that higher-rank
graph C*-algebras may be seen as special cases of our construction.
  See \cite{\KP} for definitions and a detailed treatment of
higher-rank graph C*-algebras.

Before we embark on the study of $k$-graphs from the point of view of
semigroupoids let us propose a generalization of the notion of
higher-rank graphs to semigroupoids which are not necessarily
categories.  We will not draw any conclusions based on this notion,
limiting ourselves to note that it is a natural extension of Kumjian
and Pask's interesting idea.

\definition
  \label HHSgpd
  Let $k$ be a natural number.  A rank $k$ semigroupoid, or a
$k$-semigroupoid, is a pair $(\S,d)$, where $\S$ is a semigroupoid and
  $$
  d:\S\to\N^k,
  $$
  is a function such that
  \izitem
  \zitem for every $(f,g)\in\St$, one has that $d(fg)=d(f)+d(g)$,
  \zitem if $f\in\S$, and $n,m\in\N^k$ are such that $d(f)=n+m$, there
exists a unique pair $(g,h)\in\St$ such that $d(g)=n$, $d(h)=m$, and
$gh=f$.

\bigskip
For  example,  the Markov semigroupoid is a $1$-semigroupoid, if
equipped with the word length function. 

Let $(\S,d)$ be a $k$-graph.  In particular $\S$ is a category and
hence a semigroupoid.  Under suitable hypothesis we shall now prove
that the C*-algebra of the subjacent semigroupoid is isomorphic to the
C*-algebra of $\S$, as defined by Kumjian and Pask in
\scite{\KP}{1.5}.  In particular it will follow that the dimension
function $d$ is superfluous for the definition of the corresponding
C*-algebra.

As before, if $v\in\obj{\S}$ we will denote by $\D v$ the set of
elements $f\in\S$ for which $\r(f)=v$.
  For every $n\in\N^k$  we will moreover let 
  $$
  \D v_n = \{f\in\D v: d(f) = n\}.
  $$
  We should observe that $\D v_n$ is denoted $\S^n(v)$ in \cite{\KP}.

According to \scite{\KP}{1.4}, $\S$ is said to \stress{have no
sources} if $\D v_n$ is never empty.  In case $\D v_n$ is finite for
every $v$ and $n$ one says that $\S$ is \stress{row-finite}.

Notice that the absence of sources is a much more
stringent condition than to require that $\S$ has no springs, according
to Definition \lcite{\DefineSpring}.  In fact, since $\S$ is a
category, and hence $\d(f)\in\D f$, for every $f\in\S$, we see that
$\D f\neq\emptyset$, and hence higher-rank graphs automatically have
no springs!

Below we will work under the standing hypotheses of \cite{\KP}, but we
note that our construction is meaningful regardless of these
requirements, so it would be interesting to compare our construction
with \cite{\Muhly} where these hypotheses are not required.  This
said, we suppose throughout that $\S$ is a $k$-graph for which
  $$
  0<|\D v_n|<\infty
  \for v\in\obj \S\for n\in\N^k.
  \reqmark RFNS
  $$

\state Lemma
  \label SNVisPartition
  For every object $v$ of $\S$ and every $n\in\N^k$ one has that
 $\D v_n$ is partition of $\D v$.

  \proof
  Suppose that $f,g\in\D v_n$ are such that $f\its g$. So there
are $p,q\in\S$ such that $fp=gq$.  Since $d(f)=n=d(g)$ we have
that $f=g$, by the uniqueness of the factorization.  This shows that
the elements of $\D v_n$ are pairwise disjoint.  

In order to show that $\D v_n$ is a covering of $\D v$, let
$g\in\D v$.  By \lcite{\RFNS} pick any $h\in\D{\d(g)}_n$.  Since
  $$
  d(gh) = d(g) + d(h) = d(g) + n = n + d(g),
  $$
  we may write $gh = fk$, with $d(f)=n$, and $d(k)=d(g)$.  It follows
that $f\in\D v_n$ and $g\its f$.  
  \proofend

For the convenience of the reader we now reproduce the definition of
the C*-algebra of a $k$-graph from \scite{\KP}{1.5}.

\definition 
  \label DefineCStarAlgkGraph
  Given a $k$-graph $\S$ satisfying \lcite{\RFNS},
  the C*-algebra of $\S$,  denoted by $\HHGA$, is defined to be
the universal
  C*-algebra generated by a family $\{ \rep f : f \in \S \}$ of
partial isometries satisfying:
  \izitem
  \zitem $\{ \rep v : v \in {\rm obj}(\S)\}$ is a family of
mutually orthogonal projections,
  \zitem $\rep {f g} = \rep f \rep g$ for all $f , g \in \S$ such
that $\d ( f ) = \r ( g )$,
  \zitem $\rep f^* \rep f = \rep {\d( f )}$ for all $f \in \S$,
  \zitem for every object $v$ and every $n \in {\bf N}^k$ one has 
$\displaystyle \rep v = \sum_{f \in \D v_n} \rep f \rep f^*$.

The following is certainly well known to specialists in higher-rank
graph C*-algebras:

\state  Proposition 
  \label FinalProjComute
  For every $f$ and $g$ in $\S$ one has that 
  \izitem 
  \zitem if $f\disj g$ then  $\rep f\rep f^* \perp \rep g\rep g^*$,
  \zitem $\rep f\rep f^*$ commutes with $\rep g\rep g^*$.

\proof Recall from \scite{\KP}{3.1} that whenever $n\in\N^k$ is such
that $d(f),d(g)\leq n$, we have
  $$
  
  \rep f^*\rep g = 
  \sum \rep p\rep q^*,
  $$
  where the sum extends over all pair $(p,q)$ of elements in $\S$ such
that $fp=gq$, and $d(fp)=n$.
  So
  $$
  \rep f\rep f^* \rep g\rep g^*=
  \sum_{p,q} \rep f \rep p\rep q^* \rep g^* =
  \sum_{p,q} \rep {fp}\rep {gq}^*.
  $$
  Since the last expression is symmetric with respect to $f$ and $g$,
we see that (ii) is proved.  Moreover, when $f\disj g$, it is
clear that there exist no pairs $(p,q)$ for which $fp=gq$, and hence
(i) is proved as well.
  \proofend

We shall now prove that the crucial axiom
\lcite{\DefineCStarAlgkGraph.iv} generalizes to coverings:

\state Lemma
  \label PartitionsInkGraphs
  Let $v$ be an object of $\S$. If $H$ is a finite covering of\/ $\D
v$ then
  $$
  \rep v = \bigvee_{h\in H} \rep h\rep h^*.
  $$

\proof Let $n\in\N^k$ with $n\geq d(h)$, for every $h\in H$.  For all
$f\in \D v_n$ we know that there is some $h\in H$ such that $f\its h$,
so we may write $fx=hy$, for suitable $x$ and $y$.  Since $d(f)=n\geq
d(h)$, we may write $f=f_1f_2$, with $d(f_1)=d(h)$.
  Noticing that
  $$
  f_1f_2x=hy,
  $$
  we deduce from the unique factorization property that $f_1=h$, which
amounts to saying that $h\dil f$.  We claim that this implies that 
  $\rep f\rep f^* \leq  \rep h \rep h^*$.  In fact
  $$
  \rep h \rep h^*   \rep f\rep f^* =
  \rep {f_1} \rep {f_1}^*   \rep {f_1}\rep {f_2}\rep f^* =
  \rep {f_1}\rep {f_2}\rep f^* =
  \rep f\rep f^*.
  $$
  Summarizing, we have proved that for every $f\in\D v_n$, there
exists $h\in H$, such that $\rep f\rep f^* \leq \rep h \rep h^*$.
  Therefore
  $$
  \rep v \={(\DefineCStarAlgkGraph.iv)}
  \sum_{f\in\D v_n} \rep f\rep f^* =
  \bigvee_{f\in\D v_n} \rep f\rep f^* \leq 
  \bigvee_{h\in H} \rep h\rep h^* \leq \rep v.
  $$
  from where the conclusion follows.
  \proofend

\state Theorem
  \label HHGvsSGPG
  If $\S$ is a  $k$-graph satisfying \lcite{\RFNS} then
$\HHGA$
  is *-isomorphic to $\CS$.

  \proof Throughout this proof we denote the standard generators of
$\HHGA$ by $\{\repa f\}_{f\in\S}$, together with their
initial and final projections
  $\inia f$ and 
  $\fina f$, respectively.
  Meanwhile  the standard generators of $\CSU$ will be denoted by
$\{\repc f\}_{f\in\S}$, along with with their initial and final
projections
  $\inic f$ and 
  $\finc f$.
  In particular the $\repa f$ are known to satisfy
\lcite{\DefineCStarAlgkGraph.i-iv}, while the $\repc f$ are known to
give a tight representation of the semigroupoid $\S$.

Working within the semigroupoid C*-algebra $\CSU$, we begin by
arguing that the $\repc f$ also satisfy
\lcite{\DefineCStarAlgkGraph.i-iv}.  In fact
  \lcite{\DefineCStarAlgkGraph.i} follows from
  \lcite{\ElmntFactsOnRepCat.i-ii}, while
  \lcite{\DefineCStarAlgkGraph.ii} is a consequence of
\lcite{\DefineRepTwo.ii}.  With respect to
  \lcite{\DefineCStarAlgkGraph.iii} it was proved in
  \lcite{\ElmntFactsOnRepCat.iii}.  Finally
  \lcite{\DefineCStarAlgkGraph.iv} results from the combination of
\lcite{\SNVisPartition} and \lcite{\tightInCategories.i}.

Therefore, by the universal property of $\HHGA$, there exists
a *-homomorphism 
  $$
  \Phi: \HHGA  \to \CSU ,
  $$
  such that $\Phi(\repa f)=\repc f$, for every $f\in\S$.  Evidently
the range of $\Phi$ is contained in the closed *-subalgebra of $\CSU$
generated by the $\repc f$, also known as $\CS$.

We next move our focus to the higher-rank graph algebra
$\HHGA$, and prove that the correspondence
  $$
  f\in \S \mapsto \repa f\in \HHGA 
  $$
  is a tight representation.
  Skipping the obvious  
   \lcite{\DefineRepTwo.i}
  we notice that 
   \lcite{\DefineRepTwo.ii}
  follows from \lcite{\DefineCStarAlgkGraph.ii}
  when $\d(f)=\r(g)$.  On the other hand, if $\d(f)\neq\r(g)$ we have
  $$
  \repa f \repa g =
  \repa f \repa{\d(f)} \repa{\r(g)} \repa g = 0,
  $$
  by \lcite{\DefineCStarAlgkGraph.i}.

  We next claim that the initial and final projections of the $\repa f$
commute among themselves.  That two initial projections commute
follows from \lcite{\DefineCStarAlgkGraph.i and iii}.  Speaking of the
commutativity between an 
initial projection $\inia f$ and a final projection $\fina g$, we have
that $\inia f = \repa {\d(f)}$, by \lcite{\DefineCStarAlgkGraph.iii}
and $\fin g\leq \rep {\r(g)}$ by \lcite{\DefineCStarAlgkGraph.iv}.   So
either  
  $\fina g \leq \inia f$, if $\d(f)=\r(g)$, or 
  $\fina g \perp \inia f$, if $\d(f)\neq\r(g)$, by
\lcite{\DefineCStarAlgkGraph.i}.   In any case it is clear that $\inia
f$ and $\fina g$ commute.
  That two  final projections commute is precisely the content of
\lcite{\FinalProjComute.ii}.  

Clearly \lcite{\DefineRepTwo.iii} is granted by
\lcite{\FinalProjComute.i}.
  In order to prove \lcite{\DefineRepTwo.iv} let $f,g\in\S$ with
$\d(f)=\r(g)$.  We then have that
  $$
  \inia f \fina g =
  \repa f^*\repa f \fina g =
  \repa{\d(f)} \fina g =
  \repa{\r(g)} \repa g \repa g^* =
  \repa g \repa g^* =
  \fina g.
  $$

This shows that $\repa\null$ is a representation of $\S$ in
$\HHGA$, which we will now prove to be  tight.
  For this let
  $$
  \pi:\HHGA  \to {\cal B}(H)
  $$
  be a faithful nondegenerated representation of $\HHGA$.
Through $\pi$ we will view $\HHGA$ as a subalgebra of ${\cal
B}(H)$, and hence we may consider $\repa\null$ as a representation of $\S$
on $H$.  It is clear that $\repa\null$ is nondegenerated, according to
definition \lcite{\DefineNonDegen}.
  By \lcite{\tightInCategories.ii} it is then enough to show that for
every $v\in\obj{\S}$ and every finite partition $H$ of $v$ one has
that
  $$
  \bigvee_{h\in H} \fina h = \fina v,
  $$
  but this is precisely what was proved in
\lcite{\PartitionsInkGraphs}.

By the universal property of $\CSU$ there is a *-homomorphism
  $$
  \Psi:  \CSU  \to  B(H)
  $$
  such that $\Psi(\repc f) = \repa f$, for every $f\in\S$.
Clearly  $\Psi(\CS)\subseteq\HHGA$, so   we may then view $\Psi$ and
$\Phi$ as maps
  $$
  \Psi:  \CS \to  \HHGA
  \and   
  \Phi: \HHGA  \to \CS,
  $$
  which are obviously each others inverses.
  \proofend

  \endgroup

  \begingroup
  \bigskip\bigskip 
  \font \sc = cmcsc8 \sc
  \parskip = -1pt

  Departamento de Matem\'atica 

  Universidade Federal de Santa Catarina

  88040-900 -- Florian\'opolis -- Brasil

  \eightrm exel@\kern1pt mtm.ufsc.br 

  \endgroup
  \end